\documentclass{amsart}

\usepackage{pst-plot,pst-infixplot,pstricks,graphicx}
\usepackage{amssymb,latexsym,amsthm,amsfonts,color,fancyhdr}
\usepackage{hhline}

\newtheorem{theorem}{Theorem}[section]
\newtheorem{proposition}[theorem]{Proposition}

\theoremstyle{definition}

\theoremstyle{remark}
\newtheorem{remark}{Remark}[section]

\numberwithin{equation}{section}

\begin{document}

\newcommand{\rf}[1]{(\ref{#1})}

\def\BBox{\vrule height 0.5em width 0.6em depth 0em}

\title[Orthogonal polynomials generated by a linear structure relation]
{Orthogonal polynomials generated by a linear structure relation: Inverse problem}


\author[M. Alfaro]{M. Alfaro}
\address{Departamento de Matem\'aticas and IUMA\\
Universidad de Zaragoza\\ Spain}
\email{alfaro@unizar.es}
\author[A. Pe\~na]{A. Pe\~na}
\address{Departamento de Matem\'aticas and IUMA\\
Universidad de Zaragoza\\ Spain}
\email{anap@unizar.es}
\author{J. Petronilho}
\address{CMUC and Mathematics Department\\
Universidade de Coimbra\\ Portugal}
\email{josep@mat.uc.pt}
\author[M. L. Rezola]{M. L. Rezola}
\address{Departamento de Matem\'aticas and IUMA\\
Universidad de Zaragoza\\ Spain}
\email{rezola@unizar.es}



\begin{abstract}
Let $(P_n)_n$ and $(Q_n)_n$ be two sequences of monic polynomials
linked by a type structure relation such as
$$
Q_{n}(x)+r_nQ_{n-1}(x)=P_{n}(x)+s_nP_{n-1}(x)+t_nP_{n-2}(x)\; ,
$$
where $(r_n)_n$, $(s_n)_n$ and $(t_n)_n$ are sequences of complex numbers.

Firstly, we state necessary and sufficient conditions on the parameters such that the above relation becomes non-degenerate when both sequences $(P_n)_n$ and $(Q_n)_n$ are orthogonal with respect to regular moment linear functionals ${\bf u}$ and ${\bf v}$, respectively.

Secondly, assuming that the above relation is non-degenerate and $(P_n)_n$ is an orthogonal sequence, we obtain a characterization for the orthogonality of the sequence $(Q_n)_n$ in terms of the coefficients of the polynomials $\Phi$ and $\Psi$ which appear in the rational transformation (in the distributional sense) $\Phi {\bf u}=\Psi {\bf v}\; .$

Some illustrative examples of the developed theory are presented.
\end{abstract}

\maketitle

\section{Introduction}\label{Section-Introduction}

The analysis of $M-N$ type linear structure relations involving two monic
orthogonal polynomial sequences (MOPS), $(P_n)_n$ and $(Q_n)_n$, such as
$$
Q_n(x)+\sum_{i=1}^{M-1}r_{i,n}Q_{n-i}(x)=P_n(x)+\sum_{i=1}^{N-1}s_{i,n}P_{n-i}(x)\; ,\quad n\geq0\; ,
$$
where $M$ and $N$ are fixed positive integer numbers, and
$(r_{i,n})_n$ and $(s_{i,n})_n$ are sequences of complex numbers
(and empty sum equals zero),
has been a subject of research interest in the last decades,
both from the algebraic and the analytical point of view.
For historical references, as well as a description of several aspects
focused on the interest and importance of the study of
structure relations involving linear combinations of two MOPSs,
we refer the introductory sections in the recent works \cite{Alfaro3,APRM3}
by F. Marcell\'an and three of the authors of this article,
as well as the references therein.
Such a study is also of interest in the framework
of the theory of Sobolev orthogonal polynomials, in particular
in connection with the notion of coherent pair of measures and its generalizations,
where linear structure relations involving derivatives of
at least one of the families $(P_n)_n$ and $(Q_n)_n$ appear
(see e.g. \cite{Iserles,PacoPetro,MarcioPetro}).

It is known \cite{Petronilho} that up to some natural conditions
(avoiding degenerate cases) the above $M-N$ type structure relation leads to a
rational transformation
$$
\Phi{\bf u}=\Psi{\bf v}
$$
between the regular (or quasi--definite) moment linear functionals ${\bf u}$ and ${\bf v}$
with respect to which the sequences
$(P_n)_n$ and $(Q_n)_n$ are orthogonal (respectively),
where $\Phi$ and $\Psi$ are polynomials of (exact) degrees $M-1$ and $N-1$, respectively.
As usual, $\langle{\bf w},q\rangle$ means the action of the functional ${\bf w}$
over the polynomial $q$ and the left product of a functional ${\bf w}$ (defined in the space of all polynomials) by a
polynomial $\phi$ is defined in the distributional sense, i.e.,
$\langle\phi{\bf w},p\rangle=\langle{\bf w},\phi p\rangle$ for any polynomial $p$.
In terms of the Stieltjes transforms associated with ${\bf u}$ and  ${\bf v}$,
the above relation between the functionals leads to a linear spectral transformation,
in the sense described and studied by A. Zhedanov \cite{Zhedanov},
and by V. Spiridonov, L. Vinet, and A. Zhedanov \cite{SVZ1}.
Moreover, P. Maroni \cite{MaroniANM1995} gave a characterization of
the relation between the MOPSs associated with two regular functionals ${\bf u}$ and  ${\bf v}$
fulfilling $\Phi{\bf u}=\Psi{\bf v}$.
In connection with the study of direct problems related to orthogonal polynomials associated with this kind of
modifications of linear functionals (rational modifications),
besides the work \cite{Alfaro2},
among others we also point out the works by
W. Gautschi \cite{Gautschi}, M. Sghaier and J. Alaya \cite{SA},
M.I. Bueno and F. Marcell\'an \cite{BuenoPaco}, and J.H. Lee and K.H. Kwon \cite{LeeKwon},
in particular, in the framework of the so-called Christoffel formula and its generalizations.

Concerning the above $M-N$ type structure relation,
most of the papers in the available literature  deal with
relations considering concrete values for $M$ and $N$, specially $M,N\in\{1,2,3\}$.
Indeed, the simplest relations of types $1-2$ and $2-1$
have been studied in \cite{PacoPetro}, the
$2-2$ type relation in \cite{Alfaro1} and \cite{Alfaro2},
and the more elaborated situation involving a $1-3$ type relation has been studied in \cite{APRM3}.
In addition, the $1-N$ type relation with constant coefficients
(i.e., each $(s_{i,n})_n$ is a constant sequence) has been analyzed in \cite{Alfaro3}.
In all these works a main problem stated and solved therein was the following inverse problem:
assuming that $(P_n)_n$ is a MOPS and $(Q_n)_n$
only a simple set of polynomials ---i.e., every $Q_n$ is a polynomial of degree $n$---,
to determine necessary and sufficient conditions such that $(Q_n)_n$ becomes also a MOPS.
The general $M-N$ type relations have been considered in \cite{Petronilho},
but the results were therein obtained assuming the orthogonality
of both sequences $(P_n)_n$ and $(Q_n)_n$, as well as  some additional
assumptions ensuring non-degenerate situations.
The analysis of the regularity conditions is usually a hard task,
since it involves solving systems of nonlinear difference equations,
and in general there are not available methods for solving them.
Therefore often the success depends on the
application of \emph{ad-hoc} methods for solving such systems.

In this contribution we focus on the analysis of the $M-N$
type relation with $M=2$ and $N=3$, that is
\begin{equation}\label{2-3rel}
Q_{n}(x)+r_nQ_{n-1}(x)=P_{n}(x)+s_nP_{n-1}(x)+t_nP_{n-2}(x)\; ,\quad n\geq0\; ,
\end{equation}
where $(r_n)_n$, $(s_n)_n$, and $(t_n)_n$ are sequences of complex numbers
with the conventions $r_0=s_0=t_0=t_1=0$.

Our aim is twofold. On the first hand, we determine whether (\ref{2-3rel}) is a degenerate or  a non-degenerate structure relation.  We say that the
$2-3$ type relation (\ref{2-3rel}) is degenerate if there exists another structure relation  of type $M-N$ linking $(P_n)_n$ and $(Q_n)_n$ with $M<2$ or $N<3$. In Theorem \ref{casos posibles}  we will see how some appropriate initial conditions, involving only the parameters
$r_1, r_2, r_3, s_1, s_2, t_2,$ and $t_3$, allow us to describe all the possible degenerate cases. Besides, only under the assumption $t_2\not=r_2(s_1-r_1)$ and $r_3t_3\neq0$ we really have a non--degenerate $2-3$ type relation. Since all the degenerate cases have been already considered
in the previous works \cite{Alfaro1,Alfaro2,APRM3,PacoPetro},
we will focus on the non-degenerate case. Some non-degenerate (\ref{2-3rel}) relations have been already considered in
\cite{Petronilho}, but as we will prove not all of them.

On the other hand, the following  so-called inverse problem is considered:
Given two sequences of monic polynomials $(P_n)_n$ and $(Q_n)_n$ such that (\ref{2-3rel}) holds,
and under the assumptions that (\ref{2-3rel}) is non-degenerate and $(P_n)_n$
is a MOPS, to find necessary and sufficient conditions
such that $(Q_n)_n$ becomes also a MOPS and, under such conditions,
to give the relation between the linear functionals
with respect to which $(P_n)_n$ and $(Q_n)_n$
are orthogonal. In this paper, we not only show such a characterization
(see Theorem \ref{caracterizacion ortogonalidad-1}), but we achieve
another more developed result which describes the orthogonality of $(Q_n)_n$ in terms of some sequences which remain constant (see Theorem \ref{caracterizacion ortogonalidad-2}). Even more, the most interesting fact is that the values of these constants are precisely the coefficients of the polynomials involved in the relation  between the linear functionals.
We want to observe that the same property was obtained in \cite[Theorem 2]{PacoPetro}, \cite[Theorem 2.2]{Alfaro1}, and \cite[Theorem 2.2]{APRM3}  for non--degenerate relations $1-2$, $2-2$, and $1-3$ respectively. Thus, for a general $M-N$ type structure relation (avoiding degenerate cases) we conjecture that a deeper solution of the inverse problem can be done in terms of the existence of certain constant sequences whose values coincide with the coefficients of the polynomials of (exact) degree $M-1$ and $N-1$ which relate both regular functionals.

The structure of the paper is the following.
In Section \ref{Section-degenerate}, the first of the above mentioned questions,
i.e., determining under which conditions
the $2-3$ type (\ref{2-3rel}) relation is either degenerate or non-degenerate is solved in
Theorem \ref{casos posibles}. Besides, a regularity characterization
of the functional $\Phi{\bf u}$ is given in Proposition \ref{caso (x-c)u regular},
filling out the non-degenerate relation (\ref{2-3rel}) studied
in \cite[Section 5]{Petronilho}. Also, an example of a non-degenerate
relation (\ref{2-3rel}) where the functional $\Phi{\bf u}$ is not regular is presented.
The announced regularity (orthogonality) conditions, i.e. the solution
of our inverse problem, will be stated in Theorems \ref{caracterizacion ortogonalidad-1} and \ref{caracterizacion ortogonalidad-2} in Section \ref{Section-Orthogonality}.
Finally, in Section \ref{Section-example} we present a computational example
illustrating the developed theory.
The reader may find the basic background on orthogonal polynomials needed in the sequel
in most of the articles appearing in the set of references, specially the monograph \cite{Chihara}
by T. S. Chihara where the general theory is presented, and the paper \cite{Maroni-algebrique}
by P. Maroni concerning some algebraic aspects of the theory.

\section{Degenerate and non-degenerate $2-3$ type relations}\label{Section-degenerate}

Let $(P_n)_n$ and $(Q_n)_n$ be two sequences of monic polynomials orthogonal
with respect to the regular functionals ${\bf u}$ and ${\bf v}$ (resp.),
normalized by $\langle {\bf u},1\rangle =1=\langle {\bf v},1 \rangle$.
Let $(\beta_n)_n$ and $(\gamma_n)_n$ be the sequences of recurrence coefficients
characterizing $(P_n)_n$, and $(\widetilde{\beta}_n)_n$ and $(\widetilde{\gamma}_n)_n$
the corresponding sequences characterizing $(Q_n)_n$.
Suppose that these families of polynomials are related by the $2-3$ type relation (\ref{2-3rel})
with the conventions $r_0=s_0=t_0=t_1=0$.
It is known \cite[Theorem 1.1]{Petronilho} that the initial conditions
$t_2\not=r_2(s_1-r_1)$ and $r_3t_3\not=0$
yield a relation between the linear functionals ${\bf u}$ and ${\bf v}$ such as
$$\Phi {\bf u}=\Psi {\bf v}\; ,$$
where $\Phi$ and $\Psi$ are polynomials of (exact) degrees $1$ and $2$,
respectively.

We will show that these initial conditions are not only
sufficient but  also necessary to have a non-degenerate relation,
that is $r_nt_n\not=0$ for all $n\ge 3$.

Firstly, we  point out that the conditions  $t_2\not=r_2(s_1-r_1)$ and $r_3\not=0$ imply that there
exists a complex number $c$ such that $\langle (x-c){\bf u}, Q_3 \rangle=0$ and therefore $\langle
(x-c){\bf u}, Q_n \rangle=0$ for all $n\ge 3$. Indeed, for an arbitrary $c\in\mathbb{C}$, we may write
$$
\begin{array}{rcl}
\langle (x-c){\bf u}, Q_3 \rangle &=&\langle {\bf u}, (x-c) (P_3+s_3P_2+t_3P_1-r_3Q_2) \rangle \\
&=&t_3\langle {\bf u}, (x-c)P_1 \rangle-r_3\langle {\bf u}, (x-c)[P_2+(s_2-r_2)P_1+t_2-r_2(s_1-r_1)]\rangle \\ &=&[t_3-r_3(s_2-r_2)]\gamma_1-r_3[t_2-r_2(s_1-r_1)](\beta_0-c)\;.
\end{array}
$$
Then there exists $c$ such that $\langle (x-c){\bf u}, Q_3 \rangle=0$, more precisely
\begin{equation}\label{const-c1}
c:=\beta_0- \frac{\gamma_1}{r_3}\frac{t_3-r_3(s_2-r_2)}{t_2-r_2(s_1-r_1)}\, .
\end{equation}
For this choice of $c$ and taking into account the $2-3$ type relation (\ref{2-3rel}) we have
$$\langle (x-c){\bf u}, Q_n \rangle=-r_n\langle (x-c){\bf u}, Q_{n-1}\rangle\;, \quad n\ge 4 \; .$$
Thus $\langle (x-c){\bf u}, Q_n \rangle=0$ for all $n\ge 3$, as we wish to prove.
Moreover, using standard results \cite{Maroni-algebrique},
we obtain the relation between the
functionals $(x-c){\bf u}$ and ${\bf v}$:
\begin{equation}\label{reluv1}
(x-c){\bf u}=\sum_{j=0}^2 \frac{\langle (x-c){\bf u}, Q_j \rangle}{\langle {\bf v}, Q^2_j
\rangle}\, Q_j\,{\bf v}\; ,
\end{equation}
being
\begin{align}\label{condiciones iniciales}
\langle (x-c){\bf u}, Q_0 \rangle&=\beta_0-c\;, \notag \\
\langle (x-c){\bf u}, Q_1 \rangle&=\gamma_1+(s_1-r_1)(\beta_0-c)\;, \\
\langle (x-c){\bf u}, Q_2 \rangle&=
(s_2-r_2)\gamma_1+(\beta_0-c)[t_2-r_2(s_1-r_1)]=\frac{\gamma_1t_3}{r_3}\;. \notag
\end{align}
Therefore if $t_2\not=r_2(s_1-r_1)$ and $r_3t_3\not=0$, we see that the relation between
the regular functionals ${\bf u}$ and ${\bf v}$ is $\,(x-c){\bf u}=h_2(x) {\bf v}\,$,
where $h_2$ is a polynomial of exact degree two.
Moreover, if $t_2\not=r_2(s_1-r_1)$, $r_3\not=0$, and $t_3=0$,
then $\langle (x-c){\bf u}, Q_2 \rangle=0$ and so
(\ref{reluv1}) reduces to
$(x-c){\bf u}=h_1{\bf v}$, with $h_1$ a polynomial of degree less than or equal to one,
so we have a degenerate case. In the next theorem we deduce all
the possible degenerate cases from some appropriate initials conditions
involving only the seven parameters
$r_1,r_2,r_3,s_1,s_2,t_2,$ and $t_3$.

\begin{theorem}\label{casos posibles}
Let $(P_n)_n$ and $(Q_n)_n$ be two MOPSs with
respect to the regular functionals ${\bf u}$ and ${\bf v}$, respectively, normalized by
$\langle {\bf u},1\rangle =1=\langle {\bf v},1 \rangle$. Assume that there exist sequences
of complex numbers $(r_n)_n$, $(s_n)_n$, and $(t_n)_n$ such that the $2-3$ type relation
(\ref{2-3rel}) holds, with $r_0=s_0=t_0=t_1=0$.
We have
\begin{itemize}
\item [{\rm (i)}] If $t_2=r_2(s_1-r_1)$ and $s_1=r_1$, then $t_n=0, \, n\ge 2$ and $s_n=r_n, \, n\ge1$.
Thus (\ref{2-3rel}) reduces to the trivial $1-1$ type relation $Q_n=P_n$, $n\geq0$.
\item [{\rm (ii)}] If $t_2=r_2(s_1-r_1)$ and $s_1\not=r_1$,
then $t_n=r_n(s_{n-1}-r_{n-1}), \, n\ge 2$ and $s_n\not=r_n, \, n\ge1$.
In this case (\ref{2-3rel}) reduces to a $1-2$ type relation:
$$Q_n=P_n+a_nP_{n-1}\; ,\; n\geq0\; ; \quad
a_n:=s_n-r_n\neq0\; ,\; n\geq1\; .$$
\item [{\rm (iii)}] If $t_2\not=r_2(s_1-r_1)$ and $r_3=0$, then
$t_n\not=r_n(s_{n-1}-r_{n-1}), \, n\ge 2$ and $r_n=0, \,n\ge3$.
In this case (\ref{2-3rel}) reduces to a $1-3$ type relation:
$$
\begin{array}{c}
Q_n=P_n+a_nP_{n-1}+b_nP_{n-2}\; ,\; n\geq0\; ; \\ [0.25em]
a_n:=s_n-r_n\; ,\; n\geq1\; ;\quad b_n:=t_n-r_n(s_{n-1}-r_{n-1})\neq0\; ,\; n\geq2\; .
\end{array}
$$
\item [{\rm (iv)}] If $t_2\not=r_2(s_1-r_1)$ and $r_3\not=0$, then $r_n\not=0, \,n\ge3$.
In addition:\smallskip
\begin{itemize}
\item[{\rm (iv-a)}] If $t_3=0$ and $t_2=s_2(s_1-r_1)$, then $t_n=0=s_n,\, n\ge3$.
In this case (\ref{2-3rel}) reduces to a $2-1$ type relation:
$$
\begin{array}{c}
Q_n+c_nQ_{n-1}=P_n\; ,\; n\geq0\; ;\quad c_n:=r_n-s_n\neq0\; ,\; n\geq1 \; .
\end{array}
$$
\item[{\rm (iv-b)}] If $t_3=0$ and $t_2\not=s_2(s_1-r_1)$, then $t_n=0,\, n\ge3$ and $s_n\not=0,\,n\ge3$.
In this case (\ref{2-3rel}) reduces to a $2-2$ type relation. More precisely:

 $\bullet$ If $s_1\neq r_1$, then (\ref{2-3rel}) becomes
$$
\begin{array}{c}
Q_n+c_nQ_{n-1}=P_n+d_nP_{n-1}\; ,\; n\geq0\; ;\; \\ [0.25em]
c_1-d_1:=r_1-s_1\; \\ [0.25em]
c_2:=r_2-t_2/(s_1-r_1) \; ,\quad d_2:=s_2-t_2/(s_1-r_1) \;  \\ [0.25em]
c_n:=r_n\; ,\; n\geq3\; ,\quad d_n:=s_n\; ,\; n\geq3\; .
\end{array}
$$
so that $\,c_nd_n\neq0\; ,\; n\geq1$.

$\bullet$ If $s_1=r_1$, then
$$
\begin{array}{c}
Q_1=P_1 ; \\ [0.25em]
Q_2+r_2Q_1=P_2+s_2P_1+t_2 ; \\ [0.25em]
Q_n+r_nQ_{n-1}=P_n+s_nP_{n-1}\; ,\; n\geq3\; . \\ [0.25em]
\end{array}
$$
\end{itemize}
\item[{\rm (v)}] If $t_2\not=r_2(s_1-r_1)$ and $r_3t_3\not=0$, then $r_nt_n\not=0, \, n\ge3$.
Thus in this case (\ref{2-3rel}) is a non-degenerate $2-3$ type relation.
\end{itemize}
\end{theorem}

\textbf{Proof.} From (\ref{2-3rel}) it follows that
\begin{equation}\label{uQn}
\begin{array}{l}
\langle {\bf u}, Q_1 \rangle =s_1-r_1\;,\quad \langle {\bf u}, Q_2 \rangle
=t_2-r_2(s_1-r_1)\;, \\ [0.25em]
\langle {\bf u}, Q_n \rangle =-r_n\,\langle {\bf u}, Q_{n-1} \rangle\;, \quad n\ge3\,.
\end{array}
\end{equation}

(i) If $t_2=r_2(s_1-r_1)$ and $s_1=r_1$, (\ref{uQn}) implies  $\langle {\bf u}, Q_n
\rangle=0, \, n\ge 1$, so ${\bf u}={\bf v}$. Thus $P_n=Q_n$ for all $n\geq0$ and the relation (\ref{2-3rel})
derives $s_n=r_n, \, n\ge1$, and $t_n=0, \, n\ge 2$.

(ii) If $t_2=r_2(s_1-r_1)$ and $s_1\not=r_1$, then from (\ref{uQn}) we have $
\langle {\bf u}, Q_1 \rangle\not=0$ and $\langle {\bf u}, Q_n \rangle=0, \, n\ge 2$. Hence,
the relation between the two functionals is ${\bf u}=h(x){\bf v}$ where $h$ is a
polynomial of degree one,
and so $Q_n(x)=P_n(x)+a_nP_{n-1}(x)$ for all $n\geq0$, with
$a_n\not=0,\,n\ge1$. Then the relation (\ref{2-3rel}) yields
$s_n=a_n+r_n, \, n\ge1$, and $t_n=r_na_{n-1}=r_n(s_{n-1}-r_{n-1}), \, n\ge2$.
Observe that we obtain a degenerate case, namely a $1-2$ type relation.

(iii) If $t_2\not=r_2(s_1-r_1)$ and $r_3=0$, from (\ref{uQn}) we deduce $\langle
{\bf u},Q_2 \rangle\not=0$ and $\langle {\bf u}, Q_n \rangle=0, \, n\ge 3$, so there exists
a polynomial $h$ of degree two such that ${\bf u}=h(x){\bf v}$. Thus,
$Q_n(x)=P_n(x)+a_nP_{n-1}(x)+b_nP_{n-2}(x)$ for all $n\geq0$, with $b_n\not=0,\,n\ge2$. Again, the
relation (\ref{2-3rel}) leads to $s_n=a_n+r_n,\,n\ge1$, $t_n=b_n+r_na_{n-1},\,
n\ge2$, and $r_nb_{n-1}=0, \,n\ge3$, so $r_n=0,\,n \ge3$, and $t_n\not=0,\,n \ge3$. Then we have
another degenerate case, namely a $1-3$ type relation.

(iv) If $t_2\not=r_2(s_1-r_1)$ and $r_3\not=0$,
from (\ref{uQn}) we see that $\langle {\bf u},Q_2\rangle\not=0$,
$\langle {\bf u},Q_3 \rangle\not=0$, and for each $ n\ge 4$ we have $\langle {\bf u},Q_n \rangle=0$
if $ r_n=0\,$. Assuming that there exists $n\geq4$ such that $r_n=0$,
let $n_0:=\min\{n \in\mathbb{N}\,|\,n\ge4,r_n=0\}$.
Then $\langle {\bf u},Q_n \rangle=0, \, n \ge n_0$ and
$\langle {\bf u},Q_n \rangle\not=0\,,\,2 \le n \le n_0-1$.
Hence, ${\bf u}=h(x){\bf v}$,
with $h$ a polynomial of degree $n_0-1$, and so
$$Q_n(x)=P_n(x)+\sum_{j=1}^{n_0-1}a_n^{(j)}P_{n-j}(x),$$ with
$a_n^{(n_0-1)}\not=0,\,n\ge n_0-1$.
Taking into account (\ref{2-3rel}) we easily see that this is not possible,
so $r_n\not=0, \, n\ge 3$. Moreover:\smallskip

(iv-a) If $t_3=0$ and $t_2=s_2(s_1-r_1)$, then (\ref{2-3rel}) implies $\langle {\bf v}, P_n
\rangle=0, \, n\ge 2$ and $\langle {\bf v}, P_1 \rangle\not =0$,
so ${\bf v}=h(x){\bf u}$ with $h$ a polynomial of degree one.
Then working in the same way as in (ii) we get $t_n=0=s_n, \, n\ge3$.
Note that, in this case (iv-a) we have another degenerate case, namely a $2-1$ type relation.

(iv-b) If $t_3=0$ and $t_2\not=s_2(s_1-r_1)$, then by (\ref{condiciones iniciales}) and (\ref{uQn}) we can obtain
$$
\langle (x-c){\bf u}, Q_1 \rangle=\gamma_1 \frac{t_2-s_2(s_1-r_1)}{t_2-r_2(s_1-r_1)}\not=0\; ,\quad
\langle (x-c){\bf u}, Q_n \rangle=0\,, \; n\ge 2 \; ,
$$
so there exists a polynomial $h$ of degree one such that $(x-c){\bf u}=h(x){\bf v}$.
Applying the auxiliary functional
$(x-c)^{n-2}{\bf u}$ to the main relation (\ref{2-3rel}) we obtain for $ n\ge3$
\begin{align*}
t_n \langle {\bf u}, P_{n-2 }^2\rangle
&=\langle (x-c)^{n-2}{\bf u},P_n+s_nP_{n-1}+t_nP_{n-2}\rangle \\
&=  \langle (x-c)^{n-2}{\bf u},Q_n+r_nQ_{n-1}\rangle \\
&= \langle {\bf v},(x-c)^{n-3}\,h(x)\, (Q_n+r_nQ_{n-1})\rangle=0\;.
\end{align*}
Then the condition $t_n=0,\,n\ge 3$, holds and (\ref{2-3rel}) becomes
$$Q_n(x)+r_nQ_{n-1}(x)=P_n(x)+s_nP_{n-1}(x), \, n\ge 3.$$
On the other hand, to analyze the parameters $s_n$,
we see that for $n\ge 3$
\begin{align*}
s_n \langle {\bf u}, P_{n-1}^2\rangle
&=\langle (x-c)^{n-1}{\bf u},P_n+s_nP_{n-1} \rangle
=\langle (x-c)^{n-1}{\bf u},Q_n+r_nQ_{n-1} \rangle \\
&=\langle {\bf v},(x-c)^{n-2}\,h(x)\,(Q_n+r_nQ_{n-1}) \rangle
=kr_n \langle {\bf v}, Q_{n-1}^2\rangle\; ,
\end{align*}
where $k$ is the leading coefficient of the polynomial $h$, and so we obtain $s_n\not=0, \, n\ge3$.
Thus another degenerate case appears, namely a $2-2$ type relation.

(v) If $t_2\not=r_2(s_1-r_1)$ and $r_3t_3\not=0$, then as we have seen
just before the statement of this theorem, there exists a constant $c$ such that
$(x-c){\bf u}=h_2(x){\bf v}$, with $h_2$ a  polynomial of degree two and so,
by (\ref{2-3rel}), we obtain
$$
\begin{array}{rcl}
t_n \langle {\bf u}, P_{n-2}^2 \rangle
&=&\langle (x-c){\bf u}, (P_n+s_nP_{n-1}+t_nP_{n-2})Q_{n-3} \rangle \\
&=&\langle {\bf v}, h_2(x)(Q_n+r_nQ_{n-1})Q_{n-3} \rangle
\,=\,k_2 r_n\,\langle {\bf v}, Q_{n-1}^2 \rangle\; , \;\; n\ge3\;,
\end{array}
$$
where $k_2$ is the leading coefficient of the polynomial $h_2$.
Now, it is enough to apply (iv)
to obtain $r_n\not=0, n\ge3$, and so also $t_n\not=0, n\ge3$.
Thus the proof is concluded.
\hfill$\Box$

\begin{remark}
Observe that (v) is the unique case where the relation (\ref{2-3rel})
between the two families  $(P_n)_n$ and $(Q_n)_n$ is a non-degenerate $2-3$ type relation.
All the degenerate cases, except the case (iv-b) with $r_1=s_1$, have already been considered in the previous works
\cite{Alfaro1,Alfaro2,APRM3,PacoPetro}. The case (iv-b) with $r_1=s_1$ can be studied and solved in a similar way as in \cite{Alfaro1}.
Then from now on we will concentrate on the analysis of the non-degenerate case.
\end{remark}

The non-degenerate $2-3$ type relations have already been considered in \cite[Section 5]{Petronilho}.
However, there some additional hypothesis about the parameters
involved in the relation (\ref{2-3rel}) were imposed,
namely $$t_n\not=r_n(s_{n-1}-r_{n-1})\;, \quad n\ge3\;.$$
In the following proposition we prove that to impose these
conditions, together with the conditions $r_3t_3\not=0$ and $t_2\not=r_2(s_1-r_1)$,
is equivalent to assume that the functional $(x-c){\bf u}$ is regular.

\begin{proposition}\label{caso (x-c)u regular}
Let $(P_n)_n$ and $(Q_n)_n$ be two MOPSs with
respect to the regular functionals ${\bf u}$ and ${\bf v}$, respectively, normalized by
$\langle {\bf u},1\rangle =1=\langle {\bf v},1 \rangle$. Assume that there exist sequences
of complex numbers $(r_n)_n$, $(s_n)_n$, and $(t_n)_n$ such that the $2-3$ type relation
(\ref{2-3rel}) holds, with $r_0=s_0=t_0=t_1=0$ and
the initial conditions $t_2\not=r_2(s_1-r_1)$ and $r_3t_3\not=0$.
Then the following statements are equivalent:
\begin{itemize}
\item [{\rm (i)}] The functional $(x-c){\bf u}$ is regular.
\item [{\rm (ii)}] $t_n\not=r_n(s_{n-1}-r_{n-1})$ for all $n\ge2$.
\end{itemize}
\end{proposition}

\textbf{Proof.} Multiplying both sides of (\ref{2-3rel}) by $P_{n-1}$ and applying ${\bf u}$, we find
\begin{equation}\label{uQnPn-1}
\langle {\bf u}, Q_nP_{n-1} \rangle=(s_n-r_n) \langle {\bf u}, P_{n-1}^2 \rangle\; ,\quad n\ge1\;.
\end{equation}
Moreover, multiplying both sides of (\ref{2-3rel}) by $P_{n-2}$, then applying ${\bf u}$
and taking into account (\ref{uQnPn-1}), we get
\begin{equation}\label{uQnPn-2}
\langle {\bf u}, Q_nP_{n-2} \rangle=[t_n-r_n(s_{n-1}-r_{n-1})] \langle {\bf u}, P_{n-1}^2 \rangle\;, \quad n\ge2\;.
\end{equation}
Thus
$$
t_n\not=r_n(s_{n-1}-r_{n-1}) \Longleftrightarrow
\langle {\bf u}, Q_nP_{n-2} \rangle \not=0\;, \quad n\ge2\;.
$$
On the other hand, it is well known that
$(x-c){\bf u}$ is a regular functional if and only if $P_n(c)\not=0$ for all $n\ge0$.
Therefore we only need to show that
$$\langle {\bf u}, Q_{n+2}P_n \rangle\not=0 \Longleftrightarrow  P_n(c)\not=0\;, \quad n\ge0\;.$$
Indeed, since $P_n(x)=\sum_{j=0}^n a_j^n (x-c)^j$ with $a_n^n=1$ and $a_0^n=P_n(c)$,
and (as we have seen just before the statement of Theorem \ref{casos posibles})
the relation between the regular functionals ${\bf u}$ and ${\bf v}$
is $(x-c){\bf u}=h_2(x) {\bf v}$, where $h_2$ is a polynomial of degree two, we obtain for all $n\ge1$
\begin{align*}
\langle {\bf u}, Q_{n+2}P_n\rangle
&=\langle (x-c){\bf u},Q_{n+2}[(x-c)^{n-1}+\sum_{j=1}^{n-1}a_j^n (x-c)^{j-1}] \rangle
+P_n(c) \langle {\bf u}, Q_{n+2}\rangle\\
&=\langle h_2(x) {\bf v},Q_{n+2}[(x-c)^{n-1}+\sum_{j=1}^{n-1}a_j^n (x-c)^{j-1}] \rangle
+P_n(c)\langle {\bf u}, Q_{n+2}\rangle \\
&= P_n(c)\langle {\bf u}, Q_{n+2}\rangle\;, \quad n\ge0\; .
\end{align*}
To conclude the proof it suffices to observe that from (v) in Theorem \ref{casos posibles} we have $r_n\not=0, \, n\ge3$, and therefore taking into account (\ref{uQn}) we obtain $\langle {\bf u}, Q_n\rangle\not=0$ for all $n\ge2$.
\hfill$\Box$

\medskip

Next, we are going to present an example
showing that a non-degenerate $2-3$ type relation (\ref{2-3rel}) may occur
even when the functional $(x-c){\bf u}$ is not regular.
This example shows that Theorem 5.1 in \cite{Petronilho} does not give a full description
of the non-degenerate $2-3$ type relations.
The full characterization of these relations (including the determination
of the orthogonality conditions)
is the main purpose of our study in the next section.

\medskip

\textbf{Example.}
In the sequel we denote by ${\bf w_2}$, ${\bf w_3}$, and ${\bf w_4}$ the regular functionals associated with the
Chebyshev polynomials of the second, third, and fourth kind,
which are represented (up to suitable normalizing constants) by the weight functions
$(1-x)^{1/2}(1+x)^{1/2}$, $(1-x)^{-1/2}(1+x)^{1/2}$, and
$(1-x)^{1/2}(1+x)^{-1/2}$, respectively.
Consider also the regular functional
$${\bf u}=-\mbox{$\frac{1}{3}$}\,x\,{\bf w_3}+\delta_1\; ,$$
where $\delta_\xi$ means the Dirac functional at a point $\xi\in\mathbb{C}$,
so that $\langle\delta_\xi,p\rangle:=p(\xi)$ for every polynomial $p$. The regularity of ${\bf u}$ has been stated in an example
presented in \cite[pp. 181-182]{Alfaro2}.

Denote by $(P_n)_n$, $(Q_n)_n$ and $(R_n)_n$ the MOPSs with respect
to ${\bf u}$, ${\bf w_4}$, and ${\bf w_2}$ (respectively).
These functionals satisfy the following relations
$${\bf w_2}=(1+x){\bf w_4}\;,\quad 3(x-1){\bf u}=x{\bf w_2}=x(1+x){\bf w_4}\;.$$

As a consequence of these relations and the results in \cite{Alfaro2}, we have
that the relation
$$Q_n(x)=R_n(x)+\lambda_n\,R_{n-1}(x)\; ,\quad n\ge0\,$$
with $ \lambda_n=\mbox{$\frac12$} \,, n\geq1\,,$
holds as well as the following $2-2$ type relation
$$
P_n(x)+a_nP_{n-1}(x)=R_n(x)+b_nR_{n-1}(x)\; ,\quad n\geq1\; ,
$$
where
$$a_{2n}=b_{2n}=-\frac{4n+1}{2(4n-1)}\;, \quad n\geq1\; ,$$
$$a_{2n+1}=\frac{4n-1}{2(4n+1)} \; , \quad b_{2n+1}=-\frac{4n+3}{2(4n+1)}\;, \quad n\geq0\; .$$

Moreover, since $x{\bf w_2}$ is not a regular functional,
then $(x-1){\bf u}$ is not a regular functional.
Nevertheless, a non-degenerate $2-3$ type relation between $(P_n)_n$ and $(Q_n)_n$ holds.
Indeed, noticing that $b_n\not=\lambda_n, \, n\ge1$,
we may define parameters $s_n$, $t_n$, and $r_n$, as
\begin{align*}
s_n&=a_n+\lambda_{n-1}\,\frac{b_n-\lambda_n}{b_{n-1}-\lambda_{n-1}}\, ,\quad n\geq2 \; ,\\
t_n&=a_{n-1}\,\lambda_{n-1}\,\frac{b_n-\lambda_n}{b_{n-1}-\lambda_{n-1}}\not=0\, ,\quad n\geq2 \; ,\\
r_n&=b_{n-1}\,\frac{b_n-\lambda_n}{b_{n-1}-\lambda_{n-1}}\not=0\, ,\quad n\geq2
\end{align*}
and then, by straightforward computations, we check that formula (\ref{2-3rel}) is satisfied for all $n\ge2$.
For $n=1$ we have $P_1(x)+a_1=R_1(x)+b_1=Q_1(x)-\lambda_1+b_1$,
hence $s_1-r_1=a_1-b_1+\lambda_1=\frac32\not=0$.
Moreover, $t_2=-\frac{1}{6}\neq-\frac{3}{2}= r_2(s_1-r_1)$.
Finally, we observe that $t_{2n+1}=r_{2n+1}(s_{2n}-r_{2n})$ for all $n\geq1$,
hence the conditions $t_n\not=r_n(s_{n-1}-r_{n-1})$ for all $n\ge3$ are not satisfied.

\section{Orthogonality characterizations}\label{Section-Orthogonality}

From now on, $(P_n)_n$ denotes  a MOPS with respect to a regular functional ${\bf u}$,
and $(\beta_n)_n$ and $(\gamma_{n})_n$ the corresponding sequences of recurrence
coefficients, so that
\begin{equation} \label{recurrenciaPn}
\begin{array}{l}
P_{n+1}(x)=(x-\beta_n)P_n(x)-\gamma_nP_{n-1}(x)\;, \quad n\ge0\;, \\
P_0(x)=1\;, \quad P_{-1}(x)=0\;,
\end{array}
\end{equation}
with $\gamma_n \not=0$ for all $n \ge1$.
In this section we give two characterizations of the orthogonality of
a sequence $(Q_n)_n$ of monic polynomials defined by a non-degenerate type relation (\ref{2-3rel}).
We already know from Theorem \ref{casos posibles}
that in order to have a non-degenerate $2-3$ type relation with
$(P_n)_n$ and $(Q_n)_n$ MOPSs, the conditions
$$
r_3t_3\neq0\; ,\quad t_2 \neq r_2(s_1-r_1)
$$
must hold, and these conditions imply $r_nt_n\neq0$ for all $n\geq3$.

The first characterization of the orthogonality of the sequence $(Q_n)_n$ is the following.

\begin{theorem}\label{caracterizacion ortogonalidad-1}
Let $( P_n )_{n}$ be a MOPS
and $(\beta_n)_n$ and $(\gamma_{n})_n$ the corresponding sequences of recurrence
coefficients.  We define recursively a sequence
$(Q_n)_{n}$ of monic polynomials by formula (\ref{2-3rel}), i.e.,
$$
Q_{n}(x)+r_nQ_{n-1}(x)=P_{n}(x)+s_nP_{n-1}(x)+t_nP_{n-2}(x)\;, \quad n \ge 0\; ,
$$
where $(r_n)_n$, $(s_n)_n$, and $(t_n)_n$ are sequences of complex numbers
fulfilling the conventions $r_0=s_0=t_0=t_1=0$, and such that
$$
t_2\neq r_2(s_1-r_1)\; ,\quad r_nt_n\not=0\; ,\quad n\ge3\; .
$$
Then $(Q_n)_{n}$ is a MOPS with recurrence coefficients
$(\tilde{\beta}_n)_{n}$ and $(\tilde{\gamma}_n)_{n}$, where
\begin{equation}\label{betatilden}
\widetilde{\beta}_n:=\beta_n+s_n-s_{n+1}-r_n+r_{n+1}\; ,\quad n\geq0
\end{equation}
\begin{equation}\label{gammatilden}
\widetilde{\gamma}_n:=\gamma_n+t_n-t_{n+1}+s_n(s_{n+1}-s_n-\beta_n+\beta_{n-1})
-r_n(r_{n+1}-r_n-\widetilde{\beta}_n+\widetilde{\beta}_{n-1})\; ,\quad n\geq1
\end{equation}
if and only if $\tilde{\gamma}_1\tilde{\gamma}_2\not=0$ and the following equations hold:
\begin{align}
\label{ci1} & b_2-d_2=a_2(s_1-r_1)\; , & \\
\label{ci2} & b_3-d_3=a_3(s_2-r_2)\; , & \\
\label{ci3} & c_3-b_3(s_1-r_1)=a_3[t_2-s_2(s_1-r_1)]\; , & \\
\label{eqn1} & b_n=a_ns_{n-1}\; ,\quad n\geq4\; , & \\
\label{eqn2} & c_n=a_nt_{n-1}\, ,\quad n\geq4\; , & \\
\label{eqn3} & d_n=a_nr_{n-1}\, ,\quad n\geq4\; , &
\end{align}
where
\begin{align}
\label{an} & a_n := \gamma_n+t_n-t_{n+1}+s_n(s_{n+1}-s_n-\beta_n+\beta_{n-1})\;
,\quad n\geq1\; , & \\
\label{bn} & b_n := s_n\gamma_{n-1}+t_n(s_{n+1}-s_n-\beta_n+\beta_{n-2})\;
,\quad n\geq2\; , & \\
\label{cn} & c_n := t_n\gamma_{n-2}\; ,\quad n\geq3\; , & \\
\label{dn} & d_n := r_n\widetilde{\gamma}_{n-1}\; ,\quad n\geq2\; . &
\end{align}
\end{theorem}

\textbf{Proof.} From the definition of $Q_n$ we get
\begin{equation}\label{Qn+1}
Q_{n+1}(x)=P_{n+1}(x)+s_{n+1}P_n(x)+t_{n+1}P_{n-1}(x)-r_{n+1}Q_n(x)\; , \quad n\geq0\;.
\end{equation}
Inserting formula (\ref{recurrenciaPn}) in (\ref{Qn+1}), applying (\ref{2-3rel}) to
$xP_n(x)$, and then substituting $xP_{n-1}(x)$ and $xP_{n-2}(x)$
using again (\ref{recurrenciaPn}), we get
\begin{align*}
Q_{n+1}(x)&=xQ_n(x)+(s_{n+1}-\beta_n-s_n)P_n(x)-r_{n+1}Q_n(x)+r_nxQ_{n-1}(x) \\
&\quad+(t_{n+1}-\gamma_n-s_n\beta_{n-1}-t_n)P_{n-1}(x)\\
&\quad-(s_n\gamma_{n-1}+t_n\beta_{n-2})P_{n-2}(x)-t_n\gamma_{n-2}P_{n-3}(x)\, ,\quad n\geq0\; ,
\end{align*}
with the usual convention that polynomials with negative index are zero.
Now, equation (\ref{2-3rel}) applied to $P_n(x)$ and the definition (\ref{betatilden}) of
$\tilde{\beta}_n$ yield
\begin{align*}
Q_{n+1}(x)&=(x-\tilde{\beta}_n) Q_n(x)+r_n(r_{n+1}-r_n-\tilde{\beta}_n) Q_{n-1}(x) \\
&\quad+\left[t_{n+1}-\gamma_n-t_n-s_n(s_{n+1}-s_n-\beta_n+\beta_{n-1})\right]P_{n-1}(x)\\
&\quad-\left[s_n\gamma_{n-1}+t_n(s_{n+1}-s_n-\beta_n+\beta_{n-2})\right]P_{n-2}(x)
-t_n\gamma_{n-2}P_{n-3}(x)\\
&\quad-r_n[Q_n(x)-xQ_{n-1}(x)]\; ,\quad n\geq0\; .
\end{align*}
So $(Q_n)_{n}$ is a MOPS if and only if
$\tilde{\gamma}_n \not=0$ for all $n\ge1$ and
\begin{equation}\label{auxiliar2-3}
\begin{array}{l}
r_n(r_{n+1}-r_n-\tilde{\beta}_n) Q_{n-1}(x) \\ [0.5em]
\quad+\left[ t_{n+1}-\gamma_n-t_n-s_n(s_{n+1}-s_n-\beta_n+\beta_{n-1}) \right] P_{n-1}(x) \\ [0.5em]
\quad-\left[s_n \gamma_{n-1} + t_n(s_{n+1}-s_n-\beta_n + \beta_{n-2})\right]  P_{n-2}(x)
- t_n \gamma_{n-2} P_{n-3}(x) \\ [0.5em]
\quad-r_n\left[Q_n(x)-xQ_{n-1}(x)\right]=- \tilde{\gamma}_n Q_{n-1}(x)\; ,\quad n\geq0\; .
\end{array}
\end{equation}
Moreover, $(\tilde{\beta}_n)_{n}$ and $(\tilde{\gamma}_n)_{n}$ are the corresponding sequences of recurrence coefficients of $(Q_n)_{n}.$

Next, we are going to see that $(Q_n)_{n}$ is a MOPS with recurrence coefficients
$(\tilde{\beta}_n)_{n}$ and $(\tilde{\gamma}_n)_{n}$ if and only if
$\tilde{\gamma}_n \not=0$ for all $n\ge1$ and the relation
\begin{equation}\label{nueva2-3}
\begin{array}{l}
\left[\tilde{\gamma}_n +r_n(r_{n+1}-r_n-\tilde{\beta}_n+\tilde{\beta}_{n-1}) \right] Q_{n-1}(x)
+r_n \tilde{\gamma}_{n-1} Q_{n-2}(x) \\ [0.5em]
\quad=\left[\gamma_n+ t_{n}-t_{n+1}+s_n(s_{n+1}-s_n-\beta_n+\beta_{n-1}) \right] P_{n-1}(x) \\ [0.5em]
\qquad+\left[s_n\gamma_{n-1}+t_n(s_{n+1}-s_n-\beta_{n}+\beta_{n-2}) \right]P_{n-2}(x)
+t_n\gamma_{n-2}P_{n-3}(x)
\end{array}
\end{equation}
holds for every $n\ge1$.

Suppose first that $(Q_n)_{n}$ is a MOPS with recurrence coefficients
$(\tilde{\beta}_n)_{n}$ and $(\tilde{\gamma}_n)_{n}$. Then
\begin{equation} \label{recurrenciaQn}
\begin{array}{l}
Q_{n+1}(x)=(x-\widetilde{\beta}_n)Q_n(x)-\widetilde{\gamma}_nQ_{n-1}(x)\;, \quad n\ge0\;, \\
Q_0(x)=1\;, \quad Q_{-1}(x)=0
\end{array}
\end{equation}
and $\widetilde{\gamma}_n \not=0$ for all $n \ge1$, and so
it is enough to substitute the expression for $Q_n(x)-xQ_{n-1}(x)$
obtained from this three-term recurrence relation
(after replacing $n$ by $n-1$) in formula (\ref{auxiliar2-3}) to obtain (\ref{nueva2-3}).

Conversely, if (\ref{nueva2-3}) is satisfied and $\tilde{\gamma}_n \not=0$ for all $n\ge1$,
then we show that the sequence $(Q_n)_{n}$ satisfies the three-term
recurrence relation (\ref{recurrenciaQn}),
that is,  $(Q_n)_{n}$ is a MOPS with recurrence coefficients
$(\tilde{\beta}_n)_{n}$ and $(\tilde{\gamma}_n)_{n}$.
Indeed, applying (\ref{recurrenciaPn}) in (\ref{nueva2-3}),
and by the definition of $\tilde{\beta}_n$, for $n\ge1$ we obtain
$$
\begin{array}{l}
r_n\left(\tilde{\beta}_{n-1}Q_{n-1}(x)+\tilde{\gamma}_{n-1}Q_{n-2}(x) \right)\\ [0.5em]
\quad=\gamma_nP_{n-1}(x)+(r_{n+1}-r_n-\tilde{\beta}_n)
\left[s_nP_{n-1}(x)+t_nP_{n-2}(x)-r_nQ_{n-1}(x)\right]\\ [0.5em]
\qquad+s_n\left[xP_{n-1}(x)-P_n(x)\right]+t_nxP_{n-2}(x)-t_{n+1}P_{n-1}(x)
-\tilde{\gamma}_nQ_{n-1}(x)\\ [0.5em]
\quad=\gamma_nP_{n-1}(x)+(r_{n+1}-r_n-\tilde{\beta}_n)Q_n(x)-(s_{n+1}-{\beta}_n)P_n(x) \\ [0.5em]
\qquad+x[s_nP_{n-1}(x)+t_nP_{n-2}(x)]-t_{n+1}P_{n-1}(x)-\tilde{\gamma}_nQ_{n-1}(x) \;,
\end{array}
$$
where the last equality follows from (\ref{2-3rel}).
Applying (\ref{2-3rel}) in $s_nP_{n-1}(x)+t_nP_{n-2}(x)$
as well as the recurrence relation for $(P_n)_n$, we get
$$
\begin{array}{l}
r_n\left(\tilde{\beta}_{n-1}Q_{n-1}(x)+\tilde{\gamma}_{n-1}Q_{n-2}(x) \right)
=r_n[xQ_{n-1}(x)-Q_n(x)]\\ [0.5em]
\qquad-P_{n+1}(x)+r_{n+1}Q_n(x)-s_{n+1}P_n(x)-t_{n+1}P_{n-1}(x)\\ [0.5em]
\qquad-\tilde{\beta}_nQ_n(x)+xQ_n(x)-\tilde{\gamma}_nQ_{n-1}(x)\; .
\end{array}
$$
Hence, by definition of $Q_{n+1}$ we have
$$
\begin{array}{l}
Q_{n+1}(x)-(x-\tilde{\beta}_n)Q_n(x)+\tilde{\gamma}_nQ_{n-1}(x) \\ [0.5em]
\qquad =-r_n[\,Q_n(x)-(x-\tilde{\beta}_{n-1})Q_{n-1}(x)
+\tilde{\gamma}_{n-1}Q_{n-2}(x) \,]\; ,\quad n\ge1\; .
\end{array}
$$
Therefore, since by (\ref{2-3rel})
$Q_1(x)=P_1(x)+s_1-r_1=x-\beta_0+s_1-r_1=x-\widetilde{\beta}_0$,
we deduce recursively
$$
\begin{array}{c}
Q_{n+1}(x)=(x-\tilde{\beta}_n)Q_n(x)-\tilde{\gamma}_nQ_{n-1}(x)\;,\quad n\ge0\; ,
\end{array}
$$
and so $(Q_n)_n$ is a MOPS  with recurrence coefficients
$(\tilde{\beta}_n)_{n}$ and $(\tilde{\gamma}_n)_{n}$.

Now, observe that from (\ref{2-3rel}) and (\ref{gammatilden}),
it follows that (\ref{nueva2-3}) is equivalent to
\begin{equation}\label{relacion1-2}
(d_n-r_{n-1}a_n)Q_{n-2}(x)=(b_n-s_{n-1}a_n)P_{n-2}(x)+(c_n-t_{n-1}a_n)P_{n-3}(x)\, ,
\end{equation}
for $n\geq2\, ,$ where $a_n$, $b_n$, $c_n$, and $d_n$ are defined by (\ref{an})--(\ref{dn}).

To conclude we show that (\ref{relacion1-2}) holds with
$\widetilde{\gamma}_n \not=0$ for all $n \ge1$ if and only if
$\widetilde{\gamma}_1\widetilde{\gamma}_2\neq0$ and formulas (\ref{ci1})--(\ref{eqn3}) hold.

Comparing coefficients in both sides of
(\ref{relacion1-2}) for $n=2$ and $n=3$, we obtain (\ref{ci1}) and (\ref{ci2})--(\ref{ci3}), respectively.
Moreover, it is easy to verify that (\ref{eqn1})--(\ref{eqn3}) hold for $n\ge4$.
Indeed, on the first hand, since by hypothesis $t_2\not=r_2(s_1-r_1)$ and $r_n\not=0$ for every $n\ge3$,
then by (\ref{uQn}) we deduce $\langle {\bf u}, Q_n \rangle\not=0$ for all $n\ge 2$.
On the other hand, applying ${\bf u}$ to both sides of (\ref{relacion1-2}) we obtain
$(d_n-r_{n-1}a_n)\langle {\bf u}, Q_{n-2} \rangle=0$ for $n \ge 4$. Thus $d_n-r_{n-1}a_n=0$ for every $n\geq4$,
and this proves (\ref{eqn3}). Therefore, taking into account (\ref{relacion1-2}) again, we immediately obtain
(\ref{eqn1}) and (\ref{eqn2}).

Conversely, notice that, from (\ref{eqn2}) and (\ref{eqn3}), we have
$$
\tilde{\gamma}_{n-1}=\frac{r_{n-1}}{r_n}\frac{t_n}{t_{n-1}}\gamma_{n-2}\; ,\quad n\ge4\; .
$$
Thus, $\tilde{\gamma}_n\not=0$ for every $n\ge3$, hence the conditions (\ref{ci1})--(\ref{eqn3})
together with $\widetilde{\gamma}_1\widetilde{\gamma}_2\neq0$ imply that
(\ref{relacion1-2}) holds and $\widetilde{\gamma}_n \not=0$ for all $n \ge1$.
This concludes the proof.
\hfill$\Box$

\bigskip
\begin{remark}
Using the same techniques, it can be proved a similar result exchanging the role of the sequences $(P_n)_n$ and $(Q_n)_n$. More precisely: {\it given two sequences of monic polynomials $(P_n)_n$ and $(Q_n)_n$ linked by a non--degenerate relation (\ref{2-3rel}), where  $( Q_n )_{n}$ is a MOPS with $(\tilde{\beta}_n)_{n}$ and $(\tilde{\gamma}_n)_{n}$ the corresponding sequences of recurrence coefficients, then $(P_n)_{n}$ is a MOPS with recurrence coefficients $(\beta_n)_n$ and $(\gamma_{n})_n$, satisfying (\ref{betatilden}) and (\ref{gammatilden}), if and only if $\gamma_1\gamma_2\not=0$ and the equations (\ref {ci1})--(\ref {eqn3}) hold.}
\end{remark}

\bigskip
Now, we show that the orthogonality of the sequence  $(Q_n)_n$ can be also characterized by the fact that there are three sequences (depending on the parameters $s_n, t_n, r_n$ and the recurrence coefficients) which remain constant.
Note that this new characterization of the orthogonality of $(Q_n)_n$ is more interesting than the previous one,
giving an implicit solution of the system of equations (\ref{ci1})--(\ref{eqn3}).

\begin{theorem}\label{caracterizacion ortogonalidad-2}
Let $(P_n)_n$ be a MOPS and $(\beta_n)_n$ and $(\gamma_{n})_n$
the corresponding sequences of recurrence coefficients.
Let $(Q_n)_n$ be a simple set of polynomials such that the structure relation
(\ref{2-3rel}) holds, i.e.,
$$
Q_{n}(x)+r_nQ_{n-1}(x)=P_{n}(x)+s_nP_{n-1}(x)+t_nP_{n-2}(x)\;, \quad n \ge 0\; ,
$$
where $(r_n)_n$, $(s_n)_n$ and $(t_n)_n$ are sequences of complex numbers
such that $$t_2 \neq r_2(s_1-r_1)\; ,\quad r_nt_n\neq0\; ,\quad n\geq3\; ,$$
with $r_0=s_0=t_0=t_1=0$.
Let $(\widetilde{\beta}_n)_n$ and $(\widetilde{\gamma}_{n})_n$ be defined by
(\ref{betatilden}) and (\ref{gammatilden}).
Then the following two statements are equivalent:

{\rm (i)} $(Q_n)_n$ is a MOPS with $(\widetilde{\beta}_n)_n$ and $(\widetilde{\gamma}_{n})_n$
the corresponding sequences of recurrence coefficients.

{\rm (ii)}  It holds $\widetilde{\gamma}_1\widetilde{\gamma}_2\neq0$ together with the initial conditions (\ref{ci1})--(\ref{ci3}), and
\begin{equation}\label{condicion arranque}
t_4\gamma_2=a_4t_3\; ,
\end{equation}
and the following three  sequences remain constant
\begin{equation}\label{relAn-thm}
A_n:=
\frac{s_{n}a_{n+1}}{t_{n+1}}-\beta_{n-1}-\beta_{n}+s_{n+1}=A\; ,\quad n\geq3\;
\end{equation}
\begin{equation}\label{relBn-thm}
\begin{array}{rcl}
B_n&:=&\displaystyle\frac{a_{n}a_{n+1}}{t_{n+1}}+(s_n-\beta_{n-1})
\left(\frac{s_na_{n+1}}{t_{n+1}}-\beta_{n}-s_n+s_{n+1}\right) \\ [1em]
&&\displaystyle+t_n-a_n-\gamma_{n-1}=B\; ,\quad n\geq3\; ,
\end{array}
\end{equation}
\begin{equation}\label{relCn-thm}
C_n:=
\widetilde{\beta}_n-r_{n+1}-\frac{\widetilde{\gamma}_n}{r_n}=C\; \quad n\geq3\; ,
\end{equation}
where $(a_n)_n$ is defined by (\ref{an}).
\end{theorem}

\textbf{Proof.} To prove this theorem, we introduce the auxiliary coefficients ${\beta}_n^*$ and ${\gamma}_n^*$, namely
\begin{equation}\label{betaestrellan}
\beta_n^*=\beta_n+s_n-s_{n+1}\; ,\quad n\geq0\; ,
\end{equation}
\begin{equation}\label{gamaestrellan}
\gamma_n^*=\gamma_n+t_n-t_{n+1}+s_n(s_{n+1}-s_n-\beta_n+\beta_{n-1})\; ,\quad n\geq1\; ,
\end{equation}
so we have $\gamma_n^*=a_n$ for all $n\ge1$.
Now, observe that the conditions (\ref{eqn1})--(\ref{eqn3}) in Theorem \ref{caracterizacion ortogonalidad-1} may be rewritten as
\begin{equation}\label{2.3*}
s_{n-1}\gamma_n^*=s_n\gamma_{n-1}+t_n(\beta_{n-2}-\beta_n^*)\, ,\quad n\geq4\;;
\end{equation}
\begin{equation}\label{2.4*}
t_{n-1}\gamma_n^*=t_n\gamma_{n-2}\, ,\quad n\geq4\; ;
\end{equation}
\begin{equation}\label{3*}
r_{n-1}\gamma_n^*=r_n\widetilde{\gamma}_{n-1}\, ,\quad n\geq4\; ;
\end{equation}
and therefore $(Q_n)_n$ is a MOPS if and only if
the condition $\widetilde{\gamma}_1\widetilde{\gamma}_2\neq0$,
the initial conditions (\ref{ci1})--(\ref{ci3}) and the above equations
(\ref{2.3*}), (\ref{2.4*}), and (\ref{3*}) hold.

First, since
\begin{equation*}
\gamma_n^*=a_n=\widetilde{\gamma}_n+r_n(r_{n+1}-r_n-\widetilde{\beta}_n+\widetilde{\beta}_{n-1})\; ,\quad n\geq1\; ,
\end{equation*}
we have
$$
r_n\widetilde{\gamma}_{n-1}=r_{n-1}a_n
=r_{n-1}[\widetilde{\gamma}_n+r_n(r_{n+1}-r_n-\widetilde{\beta}_n
+\widetilde{\beta}_{n-1})]\; ,\quad n\geq4\; ,
$$
hence, dividing the left and the right hand sides by $r_nr_{n-1}$,
we immediately deduce that (\ref{3*}) holds if and only if there exists a constant $C$ (independent of $n$) such that
(\ref{relCn-thm}) holds.

To conclude the proof we need to show that equations (\ref{2.3*}) and (\ref{2.4*}) are equivalent to (\ref{condicion arranque}), (\ref{relAn-thm}), and  (\ref{relBn-thm}). To do this, we notice that relations (\ref{2.3*}) and (\ref{2.4*}) are formally the same as (2.3) and (2.4) for $n\ge 4$ appearing in \cite{APRM3}, after replacing $\widetilde{\beta}_n$ and $\widetilde{\gamma}_n$ by $\beta_n^*$ and $\gamma_n^*$, respectively.

\begin{itemize}
\item {} We first prove that (\ref{2.3*}) and (\ref{2.4*})  $\Longrightarrow $ (\ref{condicion arranque}), (\ref{relAn-thm}), and (\ref{relBn-thm}):
\end{itemize}

Observe that (\ref{2.4*}) for $n=4$ is the condition (\ref{condicion arranque}). Moreover, from (\ref{2.3*}) and (\ref{2.4*}), by making exactly the same algebraic manipulations that have been
made in the proof of (i)$\Rightarrow$(ii) in
\cite[Theorem 2.2]{APRM3} we deduce that the analogue of relation (2.10) and (2.11) for $n\ge 4$ in \cite{APRM3} holds.
Thus, by straightforward computations, using the definition of $\beta_n^*$ and $\gamma_n^*$, we obtain (\ref{relAn-thm}) and (\ref{relBn-thm}).

\begin{itemize}
\item{} Next we show that (\ref{condicion arranque}), (\ref{relAn-thm}), and  (\ref{relBn-thm}) $\Longrightarrow$ (\ref{2.3*}) and (\ref{2.4*}):
\end{itemize}

As before, by making the same algebraic manipulations that have been
made in the proof of (ii)$\Rightarrow$(i) in \cite[Theorem 2.2]{APRM3},
we deduce that an equation analogous to (2.12) in \cite{APRM3} holds, i.e.,
\begin{equation*}
\frac{\gamma_{n+1}^*}{t_{n+1}}\left(\gamma_{n}-\frac{t_{n+1}}{t_{n+2}}\gamma_{n+2}^*\right)
=\gamma_{n-1}-\frac{t_n}{t_{n+1}}\gamma_{n+1}^*\, ,\quad n\geq3\, .
\end{equation*}
This relation ensures that $\gamma_{n+1}^*\neq0$ for all $n\geq3$, and so,
taking into account hypothesis (\ref{condicion arranque}), i.e., $\gamma_{n-1}-\frac{t_n}{t_{n+1}}\gamma_{n+1}^*=0$ for $n=3$, using induction we may conclude that
\begin{equation}\label{relBn-vuelta5}
\frac{\gamma_{n-1}}{t_{n}}
=\frac{\gamma_{n+1}^*}{t_{n+1}} \; ,\quad n\geq3\, .
\end{equation}
This proves (\ref{2.4*}).
To prove (\ref{2.3*}), notice first that by (\ref{relAn-thm}) we have
$$
\frac{s_{n-1}\gamma_{n}^*}{t_{n}}-\beta_{n-2}^*+s_{n-2}
=\frac{s_n\gamma_{n+1}^*}{t_{n+1}}-\beta_{n}^*+s_{n-1}\; ,\quad n\geq4\; .
$$
As a consequence, taking into account (\ref{relBn-vuelta5}),
$$
\frac{s_{n-1}\gamma_{n}^*}{t_{n}}
=\frac{s_n\gamma_{n-1}}{t_{n}}+\beta_{n-2}^*-s_{n-2}+s_{n-1}-\beta_{n}^*
=\frac{s_n\gamma_{n-1}}{t_{n}}+\beta_{n-2}-\beta_{n}^*\; ,\; n\geq4\,,
$$
and (\ref{2.3*}) holds.
Thus the proof is finished.
\hfill$\Box$

\bigskip
Next, to conclude this section, we will see that the constants $A, B$, and $C$ appearing in the above Theorem \ref{caracterizacion ortogonalidad-2} are, respectively, the coefficients $a$, $b$, and $c$ of the polynomials  which relate the two regular linear functionals, that is
\begin{equation}\label{uv}
\lambda(x-c){\bf u}=(x^2+ax+b){\bf v}.
\end{equation}

First of all, we observe that the values of $a$, $b$, $c$, and $\lambda$ may be computed from the following formulas:
\begin{equation*}\label{const-a}
a=-\widetilde{\beta}_0-\widetilde{\beta}_1
+\frac{\widetilde{\gamma}_2}{t_3}\frac{r_3t_2+(t_3-r_3s_2)(s_1-r_1)}{t_2-r_2(s_1-r_1)}\; ,
\end{equation*}
\begin{equation*}\label{const-b}
b=\widetilde{\beta}_0\widetilde{\beta}_1-\widetilde{\gamma}_1
-\frac{\widetilde{\beta}_0\widetilde{\gamma}_2}{t_3}\frac{r_3t_2+(t_3-r_3s_2)(s_1-r_1)}{t_2-r_2(s_1-r_1)}
+\frac{\widetilde{\gamma}_1\widetilde{\gamma}_2}{t_3}\frac{t_3-r_3(s_2-r_2)}{t_2-r_2(s_1-r_1)}\; ,
\end{equation*}
\begin{equation*}\label{const-c}
c=\beta_0-\frac{\gamma_1}{r_3}\frac{t_3-r_3(s_2-r_2)}{t_2-r_2(s_1-r_1)}\; ,\quad
\lambda=\frac{r_3}{t_3}\frac{\widetilde{\gamma}_1\widetilde{\gamma}_2}{\gamma_1}\; .
\end{equation*}
Indeed, making both sides of (\ref{uv}) acting on the polynomials $Q_0$, $Q_1$, and $Q_2$, and taking into account equations (\ref{condiciones iniciales}), we obtain the relations
\begin{equation}\label{eq1-abc}
\lambda(\beta_0-c)=\widetilde{\beta}_0^2+\widetilde{\beta}_0a+b+\widetilde{\gamma}_1\; ,
\end{equation}
\begin{equation}\label{eq3-abc}
\lambda[\gamma_1+(\beta_0-c)(s_1-r_1)]
=(\widetilde{\beta}_0+\widetilde{\beta}_1+a)\widetilde{\gamma}_1\; ,
\end{equation}
\begin{equation}\label{eq4-abc}
\lambda\{(s_2-r_2)\gamma_1+(\beta_0-c)[t_2-r_2(s_1-r_1)]\}
=\widetilde{\gamma}_1\widetilde{\gamma}_2\; .
\end{equation}
The expression for $c$ has been already determined in Section 2, see (\ref{const-c1}).
Hence we deduce successively $\lambda$, $a$, and $b$
from equations (\ref{eq4-abc}), (\ref{eq3-abc}), and (\ref{eq1-abc}), respectively. Note that these expressions can be also achieved by applying the proof of Theorem 1.1 in \cite{Petronilho} (which is constructive) to the particular $2-3$ type relation considered here.

\begin{theorem}\label{identificacion constantes}
Let $(P_n)_n$ and $(Q_n)_n$ be two MOPSs with respect to the regular functionals ${\bf u}$ and ${\bf v}$ respectively, normalized by $\langle {\bf u},1\rangle =1=\langle {\bf v},1 \rangle$. Let $(\beta_n,\gamma_{n+1})_{n\geq0}$ and $(\widetilde{\beta}_n,\widetilde{\gamma}_{n+1})_{n\geq0}$ be the corresponding sets of recurrence coefficients, respectively. Suppose that $(P_n)_n$ and $(Q_n)_n$ are linked by a non-degenerate 2-3 type relation such as (\ref{2-3rel}) and therefore the relation between the moment linear functionals is $$\lambda(x-c){\bf u}=(x^2+ax+b){\bf v}\; .$$
Then the constants $A$, $B$, and $C$ appearing  in Theorem \ref{caracterizacion ortogonalidad-2} coincide, respectively, with the constants $a$, $b$, and $c$.
\end{theorem}

\textbf{Proof.}

\medskip

$\bf (i) \quad C=c$

From (\ref{relCn-thm}),  $C=C_3$ and using the definition of $a_3$ in terms of $\widetilde{\gamma}_3$ (see (\ref{an}) and (\ref{gammatilden})), we have
$$C=C_3=\widetilde{\beta}_3-r_{4}-\frac{\widetilde{\gamma}_3}{r_3}=\widetilde{\beta}_2-r_{3}-\frac{a_3}{r_3}.$$

Therefore we want to prove that
$$C=\widetilde{\beta}_2-r_{3}-\frac{a_3}{r_3}=\beta_0-\frac{\gamma_1}{r_3}\frac{t_3-r_3(s_2-r_2)}{t_2-r_2(s_1-r_1)}=c.$$

Indeed, by Theorem \ref{caracterizacion ortogonalidad-1}, the initial conditions (\ref{ci1}), (\ref{ci2}), and (\ref{ci3}) hold.
Using these second and third initial conditions  we obtain
$$a_3[t_2-s_2(s_1-r_1)]=t_3\gamma_1-(s_1-r_1)[a_3(s_2-r_2)+r_3\tilde{\gamma}_2]$$
that is
\begin{equation}\label{auxiliar1-c}
t_3\gamma_1=a_3[t_2-r_2(s_1-r_1)]+r_3(s_1-r_1)\tilde{\gamma}_2.
\end{equation}
The first initial condition (\ref{ci1}) yields
$$a_2(s_1-r_1)=s_2\gamma_1+t_2(s_3-s_2-\beta_2+\beta_0)-r_2\tilde{\gamma}_1.$$
Handling adequately this expression we can obtain, see ((\ref{an}) and (\ref{gammatilden}))
\begin{equation}\label{auxiliar2-c}
(s_2-r_2)\gamma_1=(s_1-r_1)\tilde{\gamma}_2-(r_3-\widetilde{\beta}_2+\beta_0)[t_2-r_2(s_1-r_1)].
\end{equation}
Multiplying both sides of (\ref{auxiliar2-c}) by $r_3$, and then
subtracting the resulting equation from (\ref{auxiliar1-c}), we obtain
\begin{equation}\label{formula a3}
\gamma_1 \frac{t_3-r_3(s_2-r_2)}{t_2-r_2(s_1-r_1)}=a_3+r_3(r_3-\widetilde{\beta}_2+\beta_0),
\end{equation}
and then $c=C$.

$\bf (ii) \quad A=a$

First notice that from (\ref{relAn-thm}) and (\ref{condicion arranque}), and using the definition of $b_3$ (see (\ref{bn}))  we have
$$A=A_3=\frac{s_3\gamma_2}{t_3}-\beta_2-\beta_3+s_4=\frac{b_3}{t_3}-\beta_1-\beta_2+s_3.$$
We want to prove
$$A=\frac{b_3}{t_3}-\beta_1-\beta_2+s_3=-\widetilde{\beta}_0-\widetilde{\beta}_1
+\frac{\widetilde{\gamma}_2}{t_3}\frac{r_3t_2-(t_3-r_3s_2)(r_1-s_1)}{t_2-r_2(s_1-r_1)}=a.$$
Taking into account the expression of $\widetilde{\beta}_n$ (see (\ref{betatilden})) and formula (\ref{formula a3}) we get
\begin{equation}\label{auxiliar a}
A+\widetilde{\beta}_0+\widetilde{\beta}_1=\frac{b_3}{t_3}+r_3-\widetilde{\beta}_2+\beta_0=\frac{b_3}{t_3}+\frac{\gamma_1}{r_3}\frac{t_3-r_3(s_2-r_2)}{t_2-r_2(s_1-r_1)}-\frac{a_3}{r_3}.\\
\end{equation}
Now, applying successively the second initial condition (\ref{ci2}), the definition of $d_3$ (see (\ref{dn})), and  (\ref{auxiliar1-c}) we deduce
\begin{align*}
A+\widetilde{\beta}_0+\widetilde{\beta}_1&=\frac{a_3(s_2-r_2)+r_3\widetilde{\gamma}_2}{t_3}-\frac{a_3}{r_3}+
\frac{\gamma_1}{r_3}\frac{t_3-r_3(s_2-r_2)}{t_2-r_2(s_1-r_1)}\\
&=\frac{r_3}{t_3}\widetilde{\gamma}_2-\frac{a_3}{r_3t_3}(t_3-r_3(s_2-r_2))+\frac{\gamma_1}{r_3}\frac{t_3-r_3(s_2-r_2)}{t_2-r_2(s_1-r_1)}\\
&=\frac{t_3-r_3(s_2-r_2)}{r_3t_3}\left \{ \frac{t_3\gamma_1}{t_2-r_2(s_1-r_1)}-a_3 \right \}+\frac{r_3}{t_3}\widetilde{\gamma}_2\\
&=\frac{t_3-r_3(s_2-r_2)}{r_3t_3} \frac{r_3(s_1-r_1)\widetilde{\gamma}_2}{t_2-r_2(s_1-r_1)}
+ \frac{r_3}{t_3}\widetilde{\gamma}_2\;,
\end{align*}
hence $A=a$.

$\bf (iii)\quad B=b$

From (\ref{relBn-thm}) and (\ref{condicion arranque}) we have
$$B=B_3=\frac{a_3\gamma_2}{t_3}+t_3-a_3-\gamma_2+(s_3-\beta_2)(\frac{s_3\gamma_2}{t_3}-\beta_3-s_3+s_4),$$
and using the definition of $\widetilde{\beta}_n$ (see (\ref{betatilden})),  and the fact already proved $A_3=a$ we obtain
$$B=\frac{\gamma_2-t_3}{t_3}a_3+(t_3-\gamma_2)+(r_3-\widetilde{\beta}_2+s_2-r_2)(a+\widetilde{\beta}_2-r_3+r_2-s_2).$$
Hence, taking into account the expression of $b$ in terms of $\lambda$, $a$ and $c$  given by (\ref{eq1-abc}), we have to prove
$$B=-\widetilde{\gamma}_1+\lambda(\beta_0-c)-\widetilde{\beta}_0(\widetilde{\beta}_0+a).$$

Observe that by the definition of $\widetilde{\beta}_n$ and $\widetilde{\gamma}_n$
(see (\ref{betatilden}) and (\ref{gammatilden})) we can write
\begin{equation}\label{sumando1-B}
t_3-\gamma_2=t_2-r_2(s_1-r_1)-\widetilde{\gamma}_2+(s_2-r_2)(r_3-\widetilde{\beta}_2+\beta_1)
\end{equation}
and then by the equation (\ref{auxiliar1-c}) we get
\begin{equation}\label{sumando2-B}\frac{\gamma_2-t_3}{t_3}a_3=-\gamma_1+(s_1-r_1)\frac{r_3\widetilde{\gamma}_2}{t_3}+a_3\frac{\widetilde{\gamma}_2}{t_3}
-\frac{a_3(s_2-r_2)(r_3-\widetilde{\beta}_2+\beta_1)}{t_3}.
\end{equation}

Next, we analyze the last term in the expression of $B$, that is $(r_3-\widetilde{\beta}_2+s_2-r_2)(a+\widetilde{\beta}_2-r_3+r_2-s_2)$.

In the sequel of the proof, we will use the following identity
$$\beta_0-\widetilde{\beta}_0-\widetilde{\beta}_1+\beta_1+r_2-s_2=0,$$
which is a direct consequence of the definition of $\widetilde{\beta}_1$ and $\widetilde{\beta}_0$, see (\ref{betatilden}). This relation together with the first equality in (\ref{auxiliar a}) and the initial condition (\ref{ci2}) leads to
$$a+\widetilde{\beta}_2-r_3+r_2-s_2=\frac{r_3\widetilde{\gamma}_2}{t_3}+\frac{a_3(s_2-r_2)}{t_3}-\beta_1.$$
Then
\begin{equation}\label{sumando3-B}
\begin{array}l
(r_3-\widetilde{\beta}_2+s_2-r_2)(a+\widetilde{\beta}_2-r_3+r_2-s_2)\\ [0.5em]
\qquad= (r_3-\widetilde{\beta}_2+\beta_1+\beta_0-\widetilde{\beta}_0-\widetilde{\beta}_1)(a+\widetilde{\beta}_2-r_3+r_2-s_2)\\ [0.5em]
\qquad= (r_3-\widetilde{\beta}_2+\beta_1)\left(\displaystyle \frac{r_3\widetilde{\gamma}_2}{t_3}+\displaystyle\frac{a_3(s_2-r_2)}{t_3}-\beta_1\right)\\ [0.5em]
\quad\qquad +(\beta_0-\widetilde{\beta}_0-\widetilde{\beta}_1)(a+\widetilde{\beta}_2-r_3+\widetilde{\beta}_0+\widetilde{\beta}_1-\beta_1-\beta_0)\\ [0.5em]
\qquad=(r_3-\widetilde{\beta}_2+\beta_1)\left(\displaystyle\frac{r_3\widetilde{\gamma}_2}{t_3}+\displaystyle\frac{a_3(s_2-r_2)}{t_3}-\beta_1+\widetilde{\beta}_1-\beta_0+\widetilde{\beta}_0 \right)\\ [0.5em]
\quad\qquad + (\beta_0-\widetilde{\beta}_0-\widetilde{\beta}_1)(a+\widetilde{\beta}_0+\widetilde{\beta}_1-\beta_0)\; .
\end{array}
\end{equation}

Now, we add the formulas (\ref{sumando1-B}), (\ref{sumando2-B}), and (\ref{sumando3-B}). Thus, by the relation $\widetilde{\beta}_0-\beta_0+s_1-r_1=0$ and the definition of $\widetilde{\gamma}_1$ (see (\ref{gammatilden}))  we deduce

$$
\begin{array}{rcl}
B&=&-\widetilde{\gamma}_1-\widetilde{\beta}_0(a+\widetilde{\beta}_0)+(s_1-r_1)\displaystyle \frac{r_3\widetilde{\gamma}_2}{t_3}
+\displaystyle \frac{a_3\widetilde{\gamma}_2}{t_3}-\widetilde{\gamma}_2\\ [1.0em]
&&+(r_3-\widetilde{\beta}_2+\beta_1)\displaystyle \frac{r_3\widetilde{\gamma}_2}{t_3}+(\beta_0-\widetilde{\beta}_1)(a+\widetilde{\beta}_0+\widetilde{\beta}_1)\; .
\end{array}
$$
Therefore, since $C=c$ that is $\widetilde{\beta}_2-r_3-\displaystyle \frac{a_3}{r_3}=c$  we get

$$
\begin{array}{rcl}
B&=&-\widetilde{\gamma}_1-\widetilde{\beta}_0(a+\widetilde{\beta}_0)-\widetilde{\gamma}_2+\displaystyle \frac{r_3\widetilde{\gamma}_2}{t_3}(\beta_1-c+s_1-r_1)\\[1.0em]
&&+(\beta_0-\widetilde{\beta}_1)(a+\widetilde{\beta}_0+\widetilde{\beta}_1)\\[1.0em]
&=&-\widetilde{\gamma}_1-\widetilde{\beta}_0(a+\widetilde{\beta}_0)-\widetilde{\gamma}_2+\displaystyle \frac{\lambda\gamma_1}{\widetilde{\gamma}_1}
(\beta_1-c+s_1-r_1+\beta_0-\widetilde{\beta}_1)\\[1.0em]
&&+\displaystyle \frac{\lambda}{\widetilde{\gamma}_1}(\beta_0-c)\displaystyle(\beta_0-\widetilde{\beta}_1)(s_1-r_1)\; ,
\end{array}
$$
where in the second equality we have used $\lambda=\frac{r_3}{t_3}\frac{\widetilde{\gamma}_1\widetilde{\gamma}_2}{\gamma_1}$  and the relation (\ref{eq3-abc}).
Finally, it is sufficient to observe (\ref{eq4-abc}) and the definition of $\widetilde{\gamma}_1$ to conclude that $B=b$.
\hfill$\Box$

\section{Example}\label{Section-example}

A new example of a non--degenerate $2-3$ type relation (\ref{2-3rel}) is presented.

From some rational transformations of the Jacobi weight function we construct two regular functionals ${\bf u}$ and ${\bf v}$ such that  $(1+x)^2{\bf v}=(1-x){\bf u}$, in the distributional sense. Moreover we analyze when their corresponding MOPSs $(P_n)_n$ and $(Q_n)_n$  satisfy a non--degenerate $2-3$ type relation (\ref{2-3rel}), and in this case we give explicitly the parameters involved in this relation. The construction will be done in several steps.

Let ${\bf w}={\bf w}^{(\alpha, \beta)}$ be the positive definite linear functional defined by the weight function  $(1-x)^{\alpha}\,(1+x)^{\beta}\chi_{(-1,1)}$ where $\alpha>-1$ and $ \beta>-1$, and $\chi_E$ represents the characteristic function of a set $E$. Denote by $(W_n)_n$ its corresponding MOPS (we have chosen this notation instead of the classical one $(P_n^{(\alpha,\beta)})_n$ to avoid confusions with the notation $(P_n)_n$ used along all this paper). It is well known (see for instance \cite {Chihara}) that the recurrence coefficients $(\beta_n, \gamma_n)$ of $(W_n)_n$ are given by

\begin{equation*}
\beta_n = \frac{\beta^2 - \alpha^2}{(2n+\alpha+\beta)(2n+\alpha+\beta+2)}, \quad n\ge0,
\end{equation*}

\begin{equation*}
\gamma_n = \frac{4n(n+\alpha)(n+\beta)(n+\alpha+\beta)}{(2n+\alpha+\beta-1)(2n+\alpha+\beta)^2(2n+\alpha+\beta+1)}, \quad n\ge1.
\end{equation*}
Besides
$$w_0=\langle {\bf w},1 \rangle=\frac{2^{\alpha+\beta+1}\Gamma(\alpha+1)\Gamma(\beta+1)}{\Gamma(\alpha+\beta+2)},$$ and for $n\ge 1$
\begin{equation}\label{norma-Wn}
\langle {\bf w}, W_n^2 \rangle= \frac{2^{2n+\alpha+\beta+1}\Gamma(n+1)\Gamma(n+\alpha+1)\Gamma(n+\beta+1)\Gamma(n+\alpha+\beta+1)}{\Gamma(2n+\alpha+\beta+1)\Gamma(2n+\alpha+\beta+2)}.  \end{equation}
First, we consider a functional $\widetilde{\bf w}$ such that $(1-x)\widetilde{{\bf w}}={\bf w}$, that is
$$
\widetilde{{\bf w}}=(1-x)^{-1}{\bf w}+\widetilde{w}_0\delta_1.
$$
If the functional $\widetilde{{\bf w}}$ is regular and  $(\widetilde{W}_n)_n$  is the  corresponding MOPS, then there exists a sequence of complex numbers $(a_n)_n$ with $a_n\not=0$ for every  $n\ge 1$, such that
\begin{equation}\label{widetilde-Wn}
\widetilde{W}_n(x)=W_n(x)+a_nW_{n-1}(x),  \quad n\ge1.
\end{equation}
The regularity of the functional $\widetilde{{\bf w}}$ is equivalent to the parameters $a_n$ in (\ref{widetilde-Wn}) being the solution of the nonlinear difference equation
\begin{equation}\label{recurrencia-an}
\beta_n-a_{n+1}-\frac{\gamma_n}{a_n}=1, \quad n\ge1,
\end{equation}
(see Theorem 2 in \cite{PacoPetro} and its proof). Observe that there is only a free parameter, namely $a_1$.

Besides since $$\widetilde{w}_0(1-\beta_0+a_1)=w_0$$ it follows that $1-\beta_0+a_1\not=0$ and using the value of $\beta_0$ we obtain
\begin{equation}\label{condicion-a1}
2(\alpha+1)+a_1(\alpha+\beta+2)\not=0.
\end{equation}

Now, the main goal is to characterize under which conditions the functional $\widetilde{\bf w}$ is regular and to obtain the expression of the parameters $a_n$.

Through clever calculations, we can deduce that if we take $a_1\not=0$ satisfying the above condition (\ref{condicion-a1}), then for $\alpha\not=0$ the functional $\widetilde{\bf w}$ is regular if and only if

\begin{equation}\label{condicion-An}
\begin{array}{l}
A_n:=\Gamma(\alpha+1)\Gamma(\alpha+\beta+2)\Gamma(n)\Gamma(n+\beta) \\ [0.25em]
+M \Gamma(\beta+1)\Gamma(n+\alpha)\Gamma(n+\alpha+\beta)\not=0,\quad  n\ge2 ,
\end{array}
\end{equation}
where
\begin{equation*}\label{M}
M:=\frac{2\alpha\widetilde{w}_0}{w_0}-(\alpha+\beta+1)=-\frac{2(\beta+1)+a_1(\alpha+\beta+1)(\alpha+\beta+2)}{2(\alpha+1)+a_1(\alpha+\beta+2)},
\end{equation*}
and moreover, we can deduce by induction that the parameters
\begin{equation}\label{expresion-an}
a_n=\frac{-2}{(2n+\alpha+\beta)(2n+\alpha+\beta-1)}\frac{A_{n+1}}{A_n}, \quad n\ge2,
\end{equation}
are the solution of the equation (\ref{recurrencia-an}). Note that when $\alpha+\beta>-1$, (\ref{expresion-an}) is valid for $n\ge1$.

Whenever $\alpha=0$, the functional $\widetilde{\bf w}$ is regular if and only if the condition (\ref{condicion-a1}) holds and

\begin{equation}\label{condicion-tildeAn}
\begin{array}{l}
\widetilde{A}_n:=\displaystyle \frac{2 \widetilde{w}_0}{w_0}-(\beta+1) \sum_{i=1}^{n-1}\left(\frac{1}{i}+\frac{1}{\beta+i}\right)\not=0,\quad  n\ge2 .
\end{array}
\end{equation}
Besides it can be proved by induction (empty sum equals zero) that
\begin{equation}\label{expresion-tildean}
a_n=\frac{-2n(n+\beta)}{(2n+\beta)(2n+\beta-1)}\frac{\widetilde{A}_{n+1}}{\widetilde{A}_n}, \quad n\ge1.
\end{equation}

From (\ref{widetilde-Wn}) and the relation $(1-x)\widetilde{{\bf w}}={\bf w}$, we obtain
\begin{equation}\label{normawidetilde-Wn}
\langle \widetilde{\bf w}, \widetilde{W}_n^2 \rangle=-a_n \langle {\bf w},  W_{n-1}^2 \rangle,\,\quad n\ge1,
\end{equation}
and so we have an explicit expression for $\langle \widetilde{\bf w}, \widetilde{W}_n^2 \rangle$ in terms of the parameter $a_1$.

In a second step, we consider the functional ${\bf u}$ verifying $(1+x)\widetilde{{\bf w}}={\bf u}$. Then,
$$
{\bf u}=(1-x)^{-1}{\bf w}^{(\alpha, \beta+1)}+u_0\,\delta_1,
$$
where $u_0=2\widetilde{w}_0-w_0=\frac{1+\beta_0-a_1}{1-\beta_0+a_1} w_0$.
Indeed, for any polynomial $p$
$$
\begin{array}{rcl}
\langle (1+x)\widetilde{{\bf w}}, p(x) \rangle&=&\langle \widetilde{\bf w}, (1+x)p(x) \rangle \\[1.0em]
&=&\langle {\bf w}, \frac{(1+x)p(x)-2p(1)}{1-x} \rangle+2\widetilde{w}_0p(1)\\[1.0em]
&=&\langle {\bf w}, (1+x)\frac{p(x)-p(1)}{1-x} \rangle +(2\widetilde{w}_0-w_0)p(1)\\[1.0em]
&=&\langle (1-x)^{-1}\, (1+x){\bf w}^{(\alpha, \beta)}+(2\widetilde{w}_0-w_0)\delta_1, p(x) \rangle \; .
\end{array}
$$

Observe that the value of $u_0$ yields $1+\beta_0-a_1\not=0$, that is

\begin{equation}\label{condicion-a1-bis}
2(\beta+1)-a_1(\alpha+\beta+2)\not=0.
\end{equation}

Since the expression of the functional ${\bf u}$ is similar to the one of $\widetilde{{\bf w}}$, taking in mind the previous study of the regularity of the functional $\widetilde{{\bf w}}$ and exchanging $\beta$ for $\beta+1$ and $M$ for $\widetilde{M}$ where
$$\widetilde{M}:=\frac{2\alpha\,u_0}{w^{(\alpha, \beta+1)}_0}-(\alpha+\beta+2)=\frac{\alpha+\beta+2}{\beta+1}\,M,$$
we can assure that for $\alpha\not=0$ the functional ${\bf u}$ is a regular functional whenever we also impose that the conditions

\begin{equation}\label{condicion-Bn}
\begin{array}{l}
B_n:=\Gamma(\alpha+1)\,\Gamma(\alpha+\beta+3)\,\Gamma(n)\,\Gamma(n+\beta+1) \\ [0.25em]
+\,\widetilde{M}\, \Gamma(\beta+2)\,\Gamma(n+\alpha)\,\Gamma(n+\alpha+\beta+1)\not=0,\quad  n\ge1 ,
\end{array}
\end{equation}
hold. For $\alpha=0$, these conditions should be replaced by
\begin{equation}\label{condicion-tildeBn}
\begin{array}{l}
\widetilde{B}_n:=\displaystyle \frac{2\,u_0}{w^{(0, \beta+1)}_0}-(\beta+2) \sum_{i=1}^{n-1}\left(\frac{1}{i}+\frac{1}{\beta+1+i}\right)\not=0,\quad  n\ge1,
\end{array}
\end{equation}
(empty sum equals zero).

Denoting by $(P_n)_n$ the MOPS associated with this regular functional ${\bf u}$, then the following linear relation
\begin{equation}\label{widetilde2-Wn}
\widetilde{W}_n(x)=P_n(x)+b_nP_{n-1}(x), \quad n\ge1,
\end{equation}
holds, where $b_n=\langle \widetilde{\bf w}, \widetilde{W}_n^2 \rangle / \langle {\bf u},  P_{n-1}^2 \rangle$.

Furthermore, from the value of  $\langle \widetilde{\bf w}, \widetilde{W}_n^2 \rangle$ given in (\ref{normawidetilde-Wn}) and exchanging $\beta$ for $\beta+1$ and $\widetilde{w}_0$ for $u_0$, we can obtain the value of $\langle {\bf u}, P_n^2 \rangle$.
Thus, from (\ref{norma-Wn}), (\ref{expresion-an}) and (\ref{expresion-tildean}) and straightforward computations we obtain
\begin{align*}
b_1&=-a_1\,\frac{w_0}{u_0} \cr
b_n&=-a_n(n-1)(n+\alpha-1)\frac{B_{n-1}}{B_n},\quad n\ge2 ,\quad \alpha\not=0 \cr
b_n&=-a_n\frac{n-1}{n+\beta}\frac{\widetilde{B}_{n-1}}{\widetilde{B}_n},\quad n\ge2, \quad \alpha=0.
\end{align*}
Notice that we have obtained the following $2-2$ relation
\begin{equation}\label{2-2}
W_n(x)+a_nW_{n-1}(x)=P_n(x)+b_nP_{n-1}(x), \quad n\ge1,
\end{equation}
and the corresponding functionals satisfy
$$(1+x){\bf w}=(1-x){\bf u}.$$

In the last step, we consider a new functional ${\bf v}$ defined by $(1+x){\bf v}={\bf w}$, that is
$$
{\bf v}=(1+x)^{-1}{\bf w}+v_0\delta_{-1}.
$$
Again, if the functional ${\bf v}$ is regular and  $(Q_n)_n$ is the corresponding MOPS, there exists a sequence of complex numbers $(c_n)_n$ with $c_n\not=0, \quad n\ge 1,$ such that
\begin{equation}\label{Qn-Wn}
Q_n(x)=W_n(x)+c_nW_{n-1}(x),  \quad n\ge1.
\end{equation}
The functional ${\bf v}$ is regular if and only if the parameters $c_n$ satisfy
\begin{equation}\label{recurrencia-cn}
\beta_n-c_{n+1}-\frac{\gamma_n}{c_n}=-1, \quad n\ge1,
\end{equation}
(see Theorem 2 in \cite {PacoPetro}). Moreover, $$v_0(1+\beta_0-c_1)=w_0$$
hence $1+\beta_0-c_1\not=0$, and using the value of $\beta_0$ we obtain
\begin{equation}\label{condicion-c1}
2(\beta+1)-c_1(\alpha+\beta+2)\not=0
\end{equation}
Now, working in the same way as we have done before with the functional $\widetilde{{\bf w}}$, we can prove (by induction) that for $\beta\not=0$ if we take  $c_1\not=0$ satisfying the above condition (\ref{condicion-c1}) and

\begin{equation}\label{condicion-Cn}
\begin{array}{l}
C_n:=\Gamma(\beta+1)\Gamma(\alpha+\beta+2)\Gamma(n)\Gamma(n+\alpha) \\ [0.25em]
+N \Gamma(\alpha+1)\Gamma(n+\beta)\Gamma(n+\alpha+\beta)\not=0,\quad  n\ge2 ,
\end{array}
\end{equation}
where
\begin{equation*}\label{N}
N:=\frac{2\beta\,v_0}{w_0}-(\alpha+\beta+1)=-\frac{2(\alpha+1)-c_1(\alpha+\beta+1)(\alpha+\beta+2)}{2(\beta+1)-c_1(\alpha+\beta+2)},
\end{equation*}
then the functional ${\bf v}$ is regular and the parameters $c_n$ in the relation (\ref{recurrencia-cn})  are given by
\begin{equation}\label{parametros-cn}
c_n=\frac{2}{(2n+\alpha+\beta)(2n+\alpha+\beta-1)}\frac{C_{n+1}}{C_n}, \quad n\ge2 .
\end{equation}
For $\beta=0$, the functional $\widetilde{\bf v}$ is regular if and only if the condition (\ref{condicion-c1}) holds and

\begin{equation}\label{condicion-tildeCn}
\begin{array}{l}
\widetilde{C}_n:=\displaystyle \frac{2v_0}{w_0}-(\alpha+1) \sum_{i=1}^{n-1}\left(\frac{1}{i}+\frac{1}{\alpha+i}\right)\not=0,\quad  n\ge2 ,
\end{array}
\end{equation}
and besides it can be proved by induction that
\begin{equation*}\label{expresion-tildecn}
c_n=\frac{-2n(n+\alpha)}{(2n+\alpha)(2n+\alpha-1)}\frac{\widetilde{C}_{n+1}}{\widetilde{C}_n}, \quad n\ge2.
\end{equation*}

Summarizing: if we take $a_1\not=0$ and $c_1\not=0$ satisfying the conditions (\ref{condicion-a1}), (\ref{condicion-An})  or (\ref{condicion-tildeAn}) if $\alpha=0$,  (\ref{condicion-a1-bis}), (\ref{condicion-Bn})  or  (\ref{condicion-tildeBn}) if $\alpha=0$, (\ref{condicion-c1}),  and (\ref{condicion-Cn})  or (\ref{condicion-tildeCn})  if $\beta=0$, then the functionals
$${\bf u}=(1-x)^{-1}{\bf w}^{(\alpha, \beta+1)}+\frac{1+\beta_0-a_1}{1-\beta_0+a_1}w_0\,\delta_1,$$
and
$${\bf v}=(1+x)^{-1}{\bf w}^{(\alpha, \beta)}+\frac{1}{1+\beta_0-c_1}w_0\delta_{-1},$$
are regular and they are related by
$$ (1-x) {\bf u}=(1+x)^2 {\bf v} \,.$$

In general, the above relation does not imply the existence of a non--degenerate $2-3$ type relation (\ref{2-3rel}) between the sequences $(P_n)_n$ and $(Q_n)_n$ associated with the functionals ${\bf u}$ and ${\bf v}$, respectively. More precisely
we can assure that a necessary  and sufficient condition to get this type of relation is
\begin{equation}\label{annocn}
a_n\not=c_n, \quad n\ge 2.
\end{equation}
Indeed, if there exists a non--degenerate $2-3$ type relation (\ref{2-3rel}), since the functional $(1-x) {\bf u}$ is regular, from Proposition \ref{caso (x-c)u regular} we have $t_n\not=r_n(s_{n-1}-r_{n-1})$ for all $n\ge2$. Now, Theorem 5.1 in \cite{Petronilho} yields
\begin{equation*}
\langle {\bf v}, Q_n^2 \rangle \not= - \langle {\bf u}, Q_n P_{n-1} \rangle , \quad n \ge 2 .
\end{equation*}
From (\ref{widetilde2-Wn}), (\ref{2-2}), and (\ref{Qn-Wn}) we have
$$\langle {\bf u}, Q_n P_{n-1} \rangle = (b_n + c_n - a_n)\langle {\bf u}, P_{n-1}^2 \rangle,$$
and
$$\langle {\bf v}, Q_n^2 \rangle = c_n \langle {\bf w}, W_{n-1}^2 \rangle= -\frac{c_n}{a_n} \langle \widetilde{\bf w}, \widetilde{W}_n^2 \rangle = -\frac{c_n}{a_n} b_n \langle {\bf u}, P_{n-1}^2 \rangle.$$
Therefore, we get
\begin{equation*}
b_n + c_n - a_n \not= \frac{c_n}{a_n} \,b_n, \quad n \ge 2 ,
\end{equation*}
 and then (\ref{annocn}) holds.

Conversely, from (\ref{2-2}) and (\ref{Qn-Wn}) we obtain  $$Q_1(x)=P_1(x)+b_1+c_1-a_1,$$
\begin{equation}\label{Q2}
Q_2(x)+(a_2-c_2)Q_1(x)=P_2(x)+b_2P_1(x)+(a_2-c_2)c_1
\end{equation}
and straightforward computations  lead us to the following explicit non--degenerate $2-3$ type relation
$$Q_n(x)+r_nQ_{n-1}(x)=P_n(x)+s_nP_{n-1}(x)+t_nP_{n-2}(x),$$
where
\begin{align*}
r_n&=a_{n-1}\frac{a_n-c_n}{a_{n-1}-c_{n-1}} \not=0, \quad n\ge 3 , \\
s_n&=b_n+c_{n-1}\frac{a_n-c_n}{a_{n-1}-c_{n-1}},\quad  n\ge 3, \\
t_n&=b_{n-1}c_{n-1}\frac{a_n-c_n}{a_{n-1}-c_{n-1}}\not=0,\quad  n\ge 3.
\end{align*}
Observe that the condition $t_2\not=r_2(s_{1}-r_{1})$ is satisfied because $a_1 \not= b_1$.

We want to remark that in the case $a_1\not=c_1$, the above relation for $n=2$ is equivalent to the relation (\ref{Q2}).

Finally, by the sake of completeness, we show that there are a wide spectrum of free parameters $a_1$ and $c_1$ which allows us to build these examples.

For instance, taking $\alpha=\beta= 1/2, $ $a_1 \not \in (-1/2,0] \cup \{ \pm 1 \},$ and $c_1 = - a_1$, the conditions (\ref{condicion-a1}), (\ref{condicion-a1-bis}) and (\ref{condicion-c1})  are trivially satisfied, and besides, it is not difficult to verify that $A_n\not=0,\, B_n\not=0, \, C_n = A_n\not=0, \,$ for every $  n\ge1$. So the functionals
$${\bf u}={\bf w}^{(-1/2, 3/2)} - \pi \frac{1+2a_1}{1+a_1}\,\delta_1,$$
and
$${\bf v}={\bf w}^{(1/2, -1/2)} - \frac{\pi}{2} \frac{1+2a_1}{1+a_1} \delta_{-1},$$
are regular and satisfy
$$ (1-x) {\bf u}=(1+x)^2 {\bf v} \,.$$

Furthermore since $c_1 = - a_1$, from (\ref{expresion-an}) and (\ref{parametros-cn}) we deduce $c_n = - a_n, \, n\ge1.$ In this case the values of the parameters of the $2-3$ type relation are given by
\begin{align*}
r_n&=a_n, \quad n\ge 2, \\
s_n&=b_n - a_n,\quad  n\ge 2, \\
t_n&= - b_{n-1} a_n, \quad  n\ge 2,
\end{align*}
with
\begin{equation*}
a_n= - \frac{1}{2} \, \frac{1-(1+2a_1)n}{1-(1+2a_1)(n-1)}, \quad n\ge2,
\end{equation*}
and
\begin{equation*}
b_n=  -a_n \, \frac{(2n-1)(1+a_1) - (1+2a_1)(n-1)n}{(2n+1)(1+a_1) - (1+2a_1)n(n+1)}, \quad  n\ge1.
\end{equation*}

\section*{Acknowledgements}

This work has been supported by the research project MTM2012-36732-C03-02 (MEC, Spain).
J. Petronilho was also partially supported by CMUC and FCT (Portugal), through
European program COMPETE/FEDER, and the research project
PTDC/MAT/098060/2008 (FCT). The other three authors were also partially supported by DGA project E-64 (Spain).

\end{document}